\documentclass[10pt,twoside]{article}

\usepackage[centertags]{amsmath}
\usepackage{amsfonts}
\usepackage{amssymb}
\usepackage{amsthm}
\usepackage{epsfig}
\usepackage{color}
\allowdisplaybreaks[4]

\date{}

\begin{document}
\baselineskip 15pt \setcounter{page}{1}
\title{\bf \Large On the extreme order statistics for stationary Gaussian sequences subject to random missing observations
\thanks{Research supported by Innovation of Jiaxing City: a program to support the talented persons
and Project of new economy research center of Jiaxing City (No. WYZB202254).}}
\author{{\small Yuan Fang$^{1}$, Zhongquan Tan$^1$\footnote{ E-mail address:  tzq728@163.com }, Yang Yang$^{2}$}\\
\\
{\small\it
1. College of Data Science, Jiaxing University, Jiaxing 314001, PR China}\\
{\small\it
2. School of  Statistics and Data Science, Nanjing Audit University, Nanjing 211815, PR China}\\
}
 \maketitle
 \baselineskip 15pt

\begin{quote}
{\bf Abstract:}\ \  Let $\mathbf{X}=\{X_{n}\}_{n\geq 1}$ be a sequence of stationary Gaussian variables
and suppose that only some of the random variables from $\mathbf{X}$ can be observed.
In this paper, by studying the limiting properties of multidimensional exceedance point processes
for $\mathbf{X}$, we derived the joint limit distribution of extreme order statistics
for the Gaussian sequence $\mathbf{X}$ and its observed ones.
The joint limit distribution of the locations and heights of the maxima for the Gaussian sequence $\mathbf{X}$ and its observed ones are
also obtained.

{\bf Key Words:}\ \  extreme order statistics; exceedance point processes; locations; Gaussian sequences;  missing observations

{\bf AMS Classification:}\ \ Primary 60G70; secondary 60G15

\end{quote}

\section{Introduction}

Let $\mathbf{X}=\{X_{n}\}_{n\geq 1}$ be a sequence of stationary random variables and suppose that only some of the random variables from $\mathbf{X}=\{X_{n}\}_{n\geq 1}$ can be observed. Let $\boldsymbol{\varepsilon}=\{\varepsilon_{n}, n\geq1\}$ be a sequence of
Bernoulli random variables which indicates which variables are observed and is independent of $\mathbf{X}=\{X_{n}\}_{n\geq 1}$.
Let $S_{n}=\sum_{k=1}^{n}\varepsilon _{k}$ be the number of observed random variables from the set $\left\{ X_{1}, X_{2} ,\cdots , X_{n}\right\}$ and
assume that
\begin{equation*}
	\frac{S_{n} }{n} \overset{p}{\longrightarrow}\lambda,\ \ \ \ n\rightarrow\infty,\
\end{equation*}
 for some random or nonrandom variable $\lambda $.

In application fields, it is very interesting to investigate the asymptotic relation between the maxima of the sequence and the maxima of
the observed ones. More precisely, we are interested in the following asymptotic
\begin{eqnarray}
\label{1.1}
\lim_{n \to \infty} P\left\{\widetilde{M}_{n}\le u_{n}\left (x\right), M_{n} \le u_{n} \left(y\right)\right \}=\left\{\begin{matrix}
H(x,y,\lambda),\ \ \  \lambda \mbox{ is a constant};  \\
E[H(x,y,\lambda)],\ \ \ \ \ \mbox{otherwise},
\end{matrix}\right.
\end{eqnarray}
where $M_{n}=\max\left\{X_{k}, 1\leq k\leq n\right\}$, $\widetilde{M}_{n}=\max\left \{ X_{k}, \varepsilon_{k}=1, 1\leq k\leq n\right\}$ and $H$ is a nondegenerate distribution function.

If $\lambda $ is a constant, Mladenovi\'{c} and Piterbarg (2006) first proved that
(\ref{1.1}) holds with
$$H(x,y,\lambda)=G^{\lambda}\left( x \right )G^{1-\lambda } \left ( y \right )$$
 for any $x< y\in \mathbb{R}$ under a type of long dependent condition $D\left( u_{n},v_{n}\right) $ and a local dependent condition ${D}'(u_{n})$,
 where $G\left( x \right )$ is a nondegenerate distribution function.
From then on, the same problem was extended to some other cases such as the Gaussian cases in Cao and Peng (2011) and Peng et al. (2010); the autoregressive process in Glava\u{s} et al. (2017); the linear process in Glava\u{s} and Mladenovi\'{c} (2020); nonstationary random fields in Panga and Pereira (2018); and the almost sure limit theorem in Tong and Peng (2011) and Tan and Wang (2012). 

When  $\lambda $ is a random variable, Krajka (2011) showed that (\ref{1.1}) still holds with
\begin{eqnarray*}
	E\left [H\left (x,y, \lambda\right )\right] =E\left [G^{\lambda } \left ( x \right )G^{1-\lambda } \left ( y \right )  \right].
\end{eqnarray*}
Hashorva et al. (2013) investigated this problem for Gaussian sequence $\mathbf{X}=\{X_{n}\}_{n\geq 1}$ with correlation function $r_{n}=Cov\left ( X_{1},X_{n+1}\right )$ satisfying the condition
\begin{eqnarray*}
	\lim_{n \to \infty}r_{n}\ln{n} =\gamma\in[0,\infty )
\end{eqnarray*}
and proved that (\ref{1.1}) still holds with
\begin{eqnarray*}
H\left (x,y, \lambda \right ) =E\left(\int_{-\infty }^{+\infty}\exp \left(-\lambda g(x,z)-\left(1-\lambda \right)g(y,z)\right)d\Phi\left(z\right )\right),
\end{eqnarray*}
 where $g(x,z)=\exp(-x-\gamma+\sqrt{2\gamma}z)$ and $\Phi(\cdot)$ denotes the standard normal distribution function.

Some similar studies for continuous time stochastic processes can be found in
Piterbarg (2004), Tan and Tang (2014), Xu et al. (2018), Ling et al. (2018), Lu and Peng (2020) and the references therein.

In practice, we may be interested in the asymptotic relation between the extreme order statistics of some sequences and their observed ones.
It is well-known that the related limit properties of extreme order statistics can be derived directly from the associated exceedance point process
(cf. Leadbetter et al. (1983)).
Peng and Tong (2010) and Peng et al. (2019) studied the exceedance point process formed by stationary Gaussian sequences and their observed ones,
but it seems that the asymptotic relation between the extreme order statistics of the sequence and its observed
ones can not be derived from their main results, since they only considered the one-dimensional exceedance point process.
For more related studies on exceedance point process, we refer to Hsing et al. (1988), Ferreira (1999), Temido (2000), Novak (2002), Ferreira and Pereira (2012), Hashorva et al. (2014), Xiao et al. (2019), Liu and Tan (2022) and the references therein. As far as we know, the study on the exceedance point process formed by stationary Gaussian sequences and their observed ones is far from complete.

 In this paper, we continue to study the joint limiting properties of exceedance point processes formed by a stationary Gaussian sequence and its observed ones. Applying the joint limiting properties of exceedance point processes, we derived the joint asymptotic distribution between the extreme order statistics of the sequence and  its observed one. Some limit results on the locations and heights of maxima for the stationary Gaussian sequence and its observed ones are also derived.

\section{Main results}

We first give the definition of the exceedance point process formed by exceedances of  some levels by a sequence of stationary Gaussian variables $\mathbf{X}=\{X_{n}\}_{n\geq 1}$. The levels used in this paper are defined by $u_{n}(x)=a_{n}^{-1}x+b_{n}$ for any $x\in \mathbb{R}$, where
$$a_{n}=(2\ln n)^{1/2}\ \mbox{and}\ \ b_{n}=a_{n}-\frac{\ln\ln n+\ln (4\pi)}{2a_{n}}.$$

\textbf{Definition 2.1}. {\sl For a Borel set $B\subset \left ( 0,1 \right ] $, the exceedance point processes formed by $\mathbf{X}=\{X_{n}\}_{n\geq 1}$ are defined as
	\begin{eqnarray*}
	N_{n}^{x}\left ( B \right )=\sum_{\frac{j}{n}\in B }I_{\left \{ X_{j}> u_{n}(x) \right\} },
	\end{eqnarray*}
where $I_{\{\cdot\}}$ denotes the indicator function.}
	
\textbf{Definition 2.2}. {\sl  For a Borel set $B\subset \left ( 0,1 \right ] $, the exceedance point processes for observed and missed random variables from
 $\mathbf{X}=\{X_{n}\}_{n\geq 1}$ are defined as
\begin{eqnarray*}
\widetilde{N}_{n}^{x}\left ( B \right )=\sum_{\frac{j}{n}\in B }I_{\left \{ X_{j}> u_{n}(x), \varepsilon _{j}=1  \right \} } \ \ \mbox{and}\ \ \widehat{N}_{n}^{x}  \left ( B \right )=\sum_{\frac{j}{n}\in B }I_{\left \{ X_{j}> u_{n}(x)   ,\varepsilon _{j}=0  \right \} },
\end{eqnarray*}
respectively.}

Now, we state our main results.

\textbf{Theorem 2.1}. {\sl Let $\mathbf{X}=\{X_{\mathbf{n}}\}_{\mathbf{n}\geq \mathbf{1}}$ be a sequence of stationary standard (mean 0, variance 1) Gaussian variables with covariance functions $r_{n}=Cov\left ( X_{1},X_{n+1}\right )$ satisfying
\begin{eqnarray}
\label{2.1}
	\lim_{n \to \infty}r_{n}\ln_{}{n} =\gamma\in[0,\infty ).
\end{eqnarray}
Suppose that $\boldsymbol{\varepsilon}=\{\varepsilon_{n}, n\geq1\}$ is a sequence of
Bernoulli random variables  which indicates which variables in $\mathbf{X}=\{X_{\mathbf{n}}\}_{\mathbf{n}\geq \mathbf{1}}$ are observed and is independent of $\mathbf{X}=\{X_{n}\}_{n\geq 1}$.
Let $S_{n}=\sum_{k=1}^{n}\varepsilon _{k}$ be the number of observed random variables from the set $\left\{ X_{1}, X_{2} ,\cdots , X_{n}\right\}$ and satisfy
\begin{equation}
\label{2.2}
	\frac{S_{n} }{n} \overset{p}{\longrightarrow}\lambda,\ \ as\ \ n\rightarrow\infty,   \
\end{equation}
where $\lambda\in[0,1]$ a.s. is a random variable. Then, for some positive integers $k, l$,
the point process $(\widetilde{N}_{n}^{x_{1}},\ldots,\widetilde{N}_{n}^{x_{k}}, \widehat{N}_{n}^{y_{1}},\ldots,\widehat{N}_{n}^{y_{l}}) $ converges in distribution to the point process $(\widetilde{N}^{x_{1}},\ldots,\widetilde{N}^{x_{k}}, \widehat{N}^{y_{1}},\ldots ,\widehat{N}^{y_{l}})$ on $\left ( 0,1 \right ] ^{k+l} $, where $\widetilde{N}^{x}$ and $\widehat{N}^{y}$, condition on $(\lambda,\xi)$, are  Poisson processes with conditional intensity $\lambda g(x,\xi)$ and $(1-\lambda)g(y,\xi)$, respectively. Here $\xi$ is a standard Gaussian random variable independent of $\lambda$ and $g(x,\xi)=\exp(-x-\gamma+\sqrt{2\gamma}\xi)$.}

\textbf{Remark 2.1}. The point processes $\widetilde{N}_{n}^{x}$ and $\widehat{N}_{n}^{y}$ are asymptotically dependent, but, condition on $(\lambda,\xi)$, they are asymptotically independent. For the weakly dependent case, i.e., the case $\gamma=0$, the dependence is captured by the random variable $\lambda$,
while for the strongly dependent case, i.e., the case $\gamma>0$, the dependence is captured by both $\lambda$ and $\xi$.

Let $B_{i}\subset \left (0,1\right]$ be Borel sets whose boundary have zero Lebesgue measure. For application, we are interested in the joint distribution of $(\widetilde{N}^{x_{1}}(B_{1}),\ldots,\widetilde{N}^{x_{k}}(B_{k}),\widehat{N}^{y_{1}}(B_{k+1}),\ldots,\widehat{N}^{y_{l}}(B_{k+l}))$.
Since it's expression is very complicated, we only give a simple case here.

\textbf{Corollary 2.1}. {\sl Suppose that the conditions of Theorem 2.1 hold. Let $B\subset \left (0,1\right]$ be a Borel set whose boundary has zero Lebesgue measure.\\
If $x>y $, then, for $0\le k_{1}\le k_{3}$, $0\le k_{2}\le k_{4}  $, we have
\begin{eqnarray}
\label{2.3}	
&&P(\widetilde{N}^{x}(B)=k_{1},\widehat{N}^{x}(B)=k_{2},\widetilde{N}^{y}(B)=k_{3},\widehat{N}^{y}(B)=k_{4})\nonumber \\
 &&=E\int_{-\infty}^{+\infty}\bigg\{\frac{[\lambda m(B)g(x,z)]^{k_{1}}}{k_{1}!}\frac{[\lambda m(B)(g(y,z)-g(x,z))]^{k_{3}-k_{1}}}{(k_{3}-k_{1})!}\exp(-\lambda m(B)g(y,z))\nonumber \\
 &&\ \ \ \ \ \ \ \ \ \ \ \ \ \ \ \times \frac{[(1-\lambda)m(B)g(x,z)]^{k_{2}}}{k_{2}!}\frac{[(1-\lambda)m(B)(g(y,z)-g(x,z))]^{k_{4}-k_{2}}}{(k_{4}-k_{2})!}\nonumber \\
&&\ \ \ \ \ \ \ \ \ \ \ \ \ \ \ \times \exp(-(1-\lambda)m(B)g(y,z))\bigg \}\mathrm{d}\Phi(z);
\end{eqnarray}
if $x\leq y$, for $k_{1}\ge k_{3}\ge0$, $k_{2}\ge k_{4}\ge 0$, we have
\begin{eqnarray}
\label{2.4}
&&P(\widetilde{N}^{x}(B)=k_{1},\widehat{N}^{x}(B)=k_{2},\widetilde{N}^{y}(B)=k_{3},\widehat{N}^{y}(B)=k_{4})\nonumber\\
&&=E\int_{-\infty}^{+\infty}\bigg\{\frac{[\lambda m(B)g(y,z)]^{k_{3}}}{k_{3}!}
\frac{[\lambda m(B)(g(x,z)-g(y,z))]^{k_{1}-k_{3}}}{(k_{1}-k_{3})!}\exp(-\lambda m(B)g(x,z))\nonumber\\
&&\ \ \ \ \ \ \ \ \ \ \ \ \ \ \ \times\frac{[(1-\lambda)m(B)g(y,z)]^{k_{4}} }{k_{4}!} \frac{[(1-\lambda)m(B)(g(x,z)-g(y,z))]^{k_{2}-k_{4}}}{(k_{2}-k_{4})!}\nonumber\\
&&\ \ \ \ \ \ \ \ \ \ \ \ \ \ \ \times\exp(-(1-\lambda)m(B)g(x,z))\bigg \}\mathrm{d}\Phi(z),
\end{eqnarray}
where $m(B)$ denotes the Lebesgue's measure of set $B$.
}

Noting the $N_{n}^{x}=\widetilde{N}_{n}^{x}+\widehat{N}_{n}^{x}$, by the continuous mapping theorem, we have the following result which will be
used to derive the joint limit distribution for extreme order statistics for observed and missed variables.

\textbf{Corollary 2.2}. {\sl Under the conditions of Theorem 2.1,
(i). the point process $(\widetilde{N}_{n}^{x},\widehat{N}_{n}^{y}, N_{n}^{z}) $ converges in distribution to the point process $( \widetilde{N}^{x},\widehat{N}^{y}, \widetilde{N}^{z}+\widehat{N}^{z}) $ on $\left(0,1\right]^{3}$;
(ii). the point process $(\widetilde{N}_{n}^{x}, N_{n}^{y}) $ converges in distribution to the point process $( \widetilde{N}^{x}, \widetilde{N}^{y}+\widehat{N}^{y}) $ on $\left(0,1\right]^{2}$;
(iii). the point process $(\widehat{N}_{n}^{x}, N_{n}^{y}) $ converges in distribution to the point process $( \widehat{N}^{x}, \widetilde{N}^{y}+\widehat{N}^{y})$ on $\left(0,1\right]^{2}$, where $\widetilde{N}^{x}+\widehat{N}^{x}$ is a Cox process (see Page 135 Leadbetter et al. (1983) for the definition), i.e., a conditional Poisson process with conditional intensity $g(x,\xi)$.
}

Applying Corollary 2.2, we derived the following useful results for the extreme order statistics.
Let ${M}_{n}^{(m)}$ be the $m$-th  maximum of $\left\{X_{j} : j\in \mathbb{N}_{n}\right\}$ with $\mathbb{N}_{n}=\{1,2,\ldots,n\}$.
Define
$$\widetilde{M}_{n}^{(k)}=\left\{\begin{matrix}
 k\mbox{-th}\  \max\{X_{j}:\varepsilon_{j}=1,j\in \mathbb{N}_{n}\},& \mbox{if} \ \sum_{j\in \mathbb{N}_{n}}\varepsilon_{j}>k, & \\
&&\\-\infty ,\ \ \ \ \ \ \ \ \ \ \ \ \ \ \  \mbox{otherwise},
\end{matrix}\right.$$
and
$$\widehat{M}_{n}^{(l)}=\left\{\begin{matrix}
 l\mbox{-th}\  \max\{X_{j}:\varepsilon_{j}=0,j\in \mathbb{N}_{n}\},& \mbox{if} \ \sum_{j\in \mathbb{N}_{n}}\varepsilon_{j}< n-l, & \\
&&\\-\infty ,\ \ \ \ \ \ \ \ \ \ \ \ \ \ \  \mbox{otherwise}.
\end{matrix}\right.$$

\textbf{Theorem 2.2}. {\sl Suppose that the conditions of Theorem 2.1 hold.
\\(i). For any $x,y\in \mathbb{R}$,
\begin{eqnarray}
\label{2.5}
&&\lim_{n \to \infty}P\left \{ \widetilde{M}_{n}^{(k)}\le u_{n}\left ( x \right ),\widehat{M}_{n}^{(l)} \le u_{n}\left ( y \right ) \right \}\nonumber\\
&&=\sum_{s=0}^{k-1}\sum_{t=0}^{l-1}E\int_{-\infty}^{+\infty}\left \{\frac{[\lambda g(x,z)]^{s} }{s!}\exp(-\lambda g(x,z))\frac{[(1-\lambda)g(y,z)]^{t}}{t!}\exp(-(1-\lambda)g(y,z)) \right \}\mathrm{d}\Phi(z).
\end{eqnarray}
\\(ii). For $x< y\in \mathbb{R}$,
\begin{eqnarray}
\label{2.6}
&&\lim_{n\to \infty}P\left \{ \widetilde{M}_{n}^{(k)}\le u_{n}\left ( x \right ),{M}_{n}^{(m)} \le u_{n}\left ( y \right ) \right \}\nonumber\\
&&=\sum_{t=0}^{m-1}\sum_{i=0}^{m-t-1}\sum_{j=i}^{k-1}E\int_{-\infty}^{+\infty}\left \{\frac{[(1-\lambda)g(y,z)]^{t}}{t!}\exp(-(1-\lambda)g(y,z))\nonumber \right.\\
&&\ \ \ \ \ \ \ \ \ \ \ \ \ \ \ \ \ \ \ \times\left.\frac{[\lambda g(y,z)]^{i}}{i!}\frac{[\lambda(g(x,z)-g(y,z))]^{j-i} }{(j-i)!}\exp(-\lambda g(x,z))\right \}\mathrm{d}\Phi(z);
\end{eqnarray}
\\ for $x\ge y \in \mathbb{R}$,
\begin{eqnarray}
\label{2.7}
&&\lim_{n\to \infty}P\left \{ \widetilde{M}_{n}^{(k)}\le u_{n}\left ( x \right ),{M}_{n}^{(m)} \le u_{n}\left ( y \right ) \right \}\nonumber\\
&&=\sum_{t=0}^{m-1}\sum_{i=0}^{m-t-1}\sum_{j=0}^{\min\{i,k-1\}}E\int_{-\infty}^{+\infty}\left \{\frac{[(1-\lambda)g(y,z)]^{t}}{t!}\exp(-g(y,z))\nonumber \right.\\
&&\ \ \ \ \ \  \ \ \  \times\left.\frac{[\lambda g(x,z)]^{j} }{j!}\frac{[\lambda (g(y,z)-g(x,z))]^{i-j} }{(i-j)!}\right\}\mathrm{d}\Phi(z).
\end{eqnarray}
}

\textbf{Remark 2.2}. The results presented in Theorem 2.2 extend some results of Mladenovi\'{c} and Piterbarg (2006), Peng and Tong (2010), Krajka (2011) and Hashorva et al. (2013)  under the Gaussian  setting.

Let $L_{n}$, $\widetilde{L}_{n}$ and $\widehat{L}_{n}$ be the location of $M_{n}$, $\widetilde{M}_{n}$ and $\widehat{M}_{n}$, respectively. The following theorem derives the joint limit distribution for the locations and heights for the maxima from the Gaussian sequence and its observed ones, which extends Theorem 2.2 of Peng et al. (2019).

\textbf{Theorem 2.3}. {\sl Suppose that the conditions of Theorem 2.1 hold. \\
(i).  For any $0< s, t<1$ and $x,y$, we have
\begin{eqnarray}
\label{2.10}
&&\lim_{n\rightarrow\infty}P\left(\frac{1}{n}\widetilde{L}_{n}\leq s, \frac{1}{n}\widehat{L}_{n}\leq t, \widetilde{M}_{n}\leq u_{n}(x), \widehat{M}_{n}\leq u_{n}(y)\right)\nonumber\\
&&= stE\int_{-\infty}^{+\infty}\exp(-\lambda g(x,z))\exp(-(1-\lambda)g(y,z))\mathrm{d}\Phi(z).
\end{eqnarray}
(ii). For any $0< s, t<1$ and $x\leq y$, we have
\begin{eqnarray}
\label{2.11}
&&\lim_{n\rightarrow\infty}P\left(\frac{1}{n}\widetilde{L}_{n}\leq s, \frac{1}{n}L_{n}\leq t, \widetilde{M}_{n}\leq u_{n}(x), M_{n}\leq u_{n}(y)\right)\nonumber\\
&&=stE\int_{-\infty}^{+\infty}\exp(-\lambda g(x,z))\exp(-(1-\lambda)g(y,z))-\lambda \exp(-g(x,z))\mathrm{d}\Phi(z)\nonumber\\
&&\ \ \ +\min\{s,t\}E\int_{-\infty}^{+\infty}\lambda \exp(-g(x,z))\mathrm{d}\Phi(z).
\end{eqnarray}
}

\textbf{Corollary 2.3}. {\sl Suppose that the conditions of Theorem 2.1 hold. \\
(i).  For any $0< s, t<1$, we have
\begin{eqnarray}
\label{2.30}
&&\lim_{n\rightarrow\infty}P\left(\frac{1}{n}\widetilde{L}_{n}\leq s, \frac{1}{n}\widehat{L}_{n}\leq t\right)= st.
\end{eqnarray}
(ii). For any $0< s, t<1$, we have
\begin{eqnarray}
\label{2.31}
&&\lim_{n\rightarrow\infty}P\left(\frac{1}{n}\widetilde{L}_{n}\leq s, \frac{1}{n}L_{n}\leq t\right)
=st[1-E(\lambda)]+\min\{s,t\}E(\lambda).
\end{eqnarray}
}

\textbf{Remark 2.3}. Corollary 2.3 shows that $\widetilde{L}_{n}$ and $\widehat{L}_{n}$ are asymptotically independent,
but $L_{n}$ is asymptotically dependent of $\widetilde{L}_{n}$ and $\widehat{L}_{n}$.

\section{Auxiliary results}

In this section, we state and prove several lemmas which will be used in the proofs of our main results. We first borrow some notation from Krajka (2011). Let $\boldsymbol{\alpha}=\{\alpha _{n},n\ge1\}$ be nonrandom sequence taking values in $\{0,1\}$. For a random variable $\lambda$ such that $0\leq\lambda\leq1$ a.s., set
$$B_{r, t}= \begin{Bmatrix}
  & w:\lambda (w)\in \left\{\begin{matrix}\left [0,\frac{1}{2^t}   \right ],  &r=0;
  & \\ \left (\frac{r}{2^t } ,\frac{r+1}{2^t }    \right ], &0<r\leq 2^t-1
\end{matrix}\right.&
\end{Bmatrix},$$
$$B_{r, t, \boldsymbol{\alpha}, n}=\{w: \varepsilon_j(w)= \alpha_j, 1 \leq j \leq n\}\cap B_{r, t}.$$
For the arbitrary random or nonrandom sequence $\boldsymbol{\beta} =\{\beta_n,n\geq 1\}$ of 0 and 1 and subset $I\subset \mathbb{N}_{n}$, define
$$\widetilde{M}_{n}(X,I,\boldsymbol{\beta})=\left\{\begin{matrix}
  \max\{X_{j}:\beta _{j}=1,j\in I\},& \mbox{if} \ \max_{j\in I}\beta_{j}>0,& \\
&&\\-\infty ,\ \ \ \ \ \ \ \ \ \ \ \ \ \ \  \mbox{otherwise},
\end{matrix}\right.$$
and
$$\widehat{M}_{n}(X,I,\boldsymbol{\beta})=\left\{\begin{matrix}
  \max\{X_{j}:\beta_{j}=0,j\in I\},& \mbox{if} \ \min_{j\in I}\beta_{j}=0,& \\
&&\\-\infty ,\ \  \ \ \ \ \ \ \ \ \ \ \ \ \  \mbox{otherwise}.
\end{matrix}\right.$$
For any integer $k>0$, let $c_i, d_i, i=1,2,\ldots,k$ be constants such that $0< c_{1}<d_{1}\le c_{2}<\cdots \le c_{k}<d_{k}\le 1$.
Let $A_{i}=\left \{\lfloor nc_{i}\rfloor+1,\lfloor nc_{i}\rfloor+2,\cdots ,\lfloor nd_{i}\rfloor\right \}$, $1\leq i\leq k$, where $\lfloor x\rfloor $ denotes the integral part of $x$. Let also $x_{i}$ and $y_{i}$, $i=1,2,\ldots,k$, be two sequences of real numbers.

\textbf{Lemma 3.1}. {\sl Let $\{Z_{j}\}_{j\ge1} $ be a sequence of independent standard Gaussian random variables. Under the assumptions of Theorem 2.1,  we have
\begin{eqnarray}
\label{10}
&&\bigg|P\left (\bigcap_{i=1}^{k}\left \{\widetilde{M}_{n}(X,A_{i},\boldsymbol{\alpha})\le u_{n}(x_{i}),\widehat{M}_{n}(X,A_{i},\boldsymbol{\alpha})\le u_{n}(y_{i})\right \}\right )\nonumber \\
 &&\ \ \ \ \ \ \ \ \ \ -\int_{-\infty}^{+\infty}\prod_{i=1}^{k}P\left \{\widetilde{M}_{n}(Z,A_{i},\boldsymbol{\alpha})\le u_{n}(x_{i},z),\widehat{M}_{n}(Z,A_{i},\boldsymbol{\alpha})\le u_{n}(y_{i},z)  \right \}\mathrm{d}\Phi(z)\bigg|\rightarrow 0,
\end{eqnarray}
as $n\rightarrow\infty $, where $u_{n}(x,z)=(1-\rho _{n})^{-\frac{1}{2}}(u_{n}(x)-\rho _{n}^{\frac{1}{2}}z)$ with $\rho _{n}=\gamma /\ln{n}$.
}

\textbf{Proof:} Define $Y_{j}=(1-\rho _{n})^{\frac{1}{2}}Z_{j}+\rho _{n}^{\frac{1}{2}}\xi$, where $\xi$ is a standard Gaussian random variable. By Normal Comparision Lemma (see e.g., Leadbetter et al. (1983)), we have
\begin{eqnarray}
\label{11}
&&\bigg|P\left (\bigcap_{i=1}^{k}\left \{\widetilde{M}_{n}(X,A_{i},\boldsymbol{\alpha})\le u_{n}(x_{i}),\widehat{M}_{n}(X,A_{i},\boldsymbol{\alpha})\le u_{n}(y_{i})\right \}\right )\nonumber \\
&&\ \ \ \ -P\left (\bigcap_{i=1}^{k}\left \{\widetilde{M}_{n}(Y,A_{i},\boldsymbol{\alpha})\le u_{n}(x_{i}),\widehat{M}_{n}(Y,A_{i},\boldsymbol{\alpha})\le u_{n}(y_{i})\right \}\right )\bigg|\nonumber \\
&&\le Kn\sum_{i=1}^{\lfloor nd_{k}\rfloor}|r_{j}-\rho _{n}|\exp\left (-\frac{(\min\{u_{n}(x_{i}),u_{n}(y_{i}),i=1,2,\ldots,k\}) ^{2}}{1+\omega _{j}}\right ),
\end{eqnarray}
where $\omega _{j}=\max\{|r_{j}|,\rho _{n}\}$ and $K$ is some constant, depending only on $\gamma$. It is directly to check that
\begin{eqnarray}
\label{12}
&&P\left (\bigcap_{i=1}^{k}\left \{\widetilde{M}_{n}(Y,A_{i},\boldsymbol{\alpha})\le u_{n}(x_{i}),\widehat{M}_{n}(Y,A_{i},\boldsymbol{\alpha})\le u_{n}(y_{i})\right \}\right )\nonumber \\
&&=\int_{-\infty}^{+\infty}P\left (\bigcap_{i=1}^{k}\left \{\widetilde{M}_{n}(Z,A_{i},\boldsymbol{\alpha})\le u_{n}(x_{i},z),\widehat{M}_{n}(Z,A_{i},\boldsymbol{\alpha})\le u_{n}(y_{i},z)\right \}\right ) \mathrm{d}\Phi(z)\nonumber \\
&&=\int_{-\infty}^{+\infty}\prod_{i=1}^{k}P\left \{\widetilde{M}_{n}(Z,A_{i},\boldsymbol{\alpha})\le u_{n}(x_{i},z),\widehat{M}_{n}(Z,A_{i},\boldsymbol{\alpha})\le u_{n}(y_{i},z) \right \}\mathrm{d}\Phi(z).
\end{eqnarray}
Applying Lemma 6.4.1 of Leadbetter et al. (1983), \eqref{11} and \eqref{12}, we obtain the desired result.

The following lemma is from  Krajka (2011).

\textbf{Lemma 3.2}. {\sl Let $d(X,Y)$ stand for Ky Fan metric, $d(X,Y)=\inf\{\varepsilon : P[|X-Y|>\varepsilon]<\varepsilon\}$.

(a). For arbitrary positive integers $s$, $m$, we have
$$d\left ( \frac{S_{ms}-S_{m(s-1)}}{m} ,\lambda\right)\le (2s-1)\left [d\left ( \frac{S_{ms}}{ms} ,\lambda\right)+d\left ( \frac{S_{m(s-1)}}{m(s-1)} ,\lambda\right) \right].$$

(b). If $\{X_{n},n\ge 1\}$ and $\{Y_{n},n\ge 1\}$ are such that $|X_{n}-Y_{n}|<1$ a.s. then
\begin{eqnarray*}
E|X_{n}-Y_{n}|\le 2d(X_{n},Y_{n}).
\end{eqnarray*}}

\textbf{Lemma 3.3}. {\sl Under the conditions of Theorem 2.1, we have
\begin{eqnarray}
\label{13}
&&P\left (\bigcap_{i=1}^{k}\left \{\widetilde{M}_{n}(X,A_{i},\boldsymbol
{\varepsilon})\le u_{n}(x_{i}),\widehat{M}_{n}(X,A_{i},\boldsymbol
{\varepsilon})\le u_{n}(y_{i})\right \}\right)\nonumber \\
&&\longrightarrow E\int_{-\infty}^{+\infty}\prod_{i=1}^{k}\exp\left [-(d_{i}-c_{i})(\lambda g(x_{i},z)+(1-\lambda)g(y_{i},z)) \right ]\mathrm{d}\Phi(z),
\end{eqnarray}
as $n \rightarrow \infty$.}

\textbf{Proof:} Let $t$ be a fixed integer and $m_{i}=\left \lfloor\frac{nd_{i}-nc_{i}}{t}\right\rfloor, i=1,2,\ldots,k$. Put
$$K_{s,i}=\{j: \lfloor nc_{i}\rfloor +(s-1)m_{i}+1\le j\le \lfloor nc_{i}\rfloor +sm_{i}\},\ \ \ \ 1\le s\le t, $$
and
$$I_{s,i}=\{j: \lfloor nc_{i}\rfloor+(s-1)m_{i}+1\le j\le \lfloor nc_{i}\rfloor+sm_{i}-l\},\ \ J_{s,i}=K_{s,i}\setminus I_{s,i}. $$
Obviously,
\begin{eqnarray}
\label{14}
&&\bigg|P\left (\bigcap_{i=1}^{k}\left \{\widetilde{M}_{n}(X,A_{i},\boldsymbol{\varepsilon})\le u_{n}(x_{i}),\widehat{M}_{n}(X,A_{i},\boldsymbol
{\varepsilon})\le u_{n}(y_{i})\right \}\right)\nonumber \\
&&\ \ \ \ \ -E\int_{-\infty}^{+\infty}\prod_{i=1}^{k}\exp\left [-(d_{i}-c_{i})(\lambda g(x_{i},z)+(1-\lambda)g(y_{i},z)) \right ]\mathrm{d}\Phi(z)\bigg|\nonumber \\
&&\le \sum_{r=0}^{2^{t}-1}\sum_{\boldsymbol{\alpha} \in \{0,1\}^{n}}E\bigg|P\left (\bigcap_{i=1}^{k}\left \{\widetilde{M}_{n}(X,A_{i},\boldsymbol{\alpha} )\le u_{n}(x_{i}),\widehat{M}_{n}(X,A_{i},\boldsymbol{\alpha} )\le u_{n}(y_{i})\right \}\right)\nonumber \\
&&\ \ \ \ \ - \int_{-\infty}^{+\infty}\prod_{i=1}^{k}\exp\left [-(d_{i}-c_{i})(\lambda g(x_{i},z)+(1-\lambda)g(y_{i},z)) \right ]\mathrm{d}\Phi(z)\bigg|I_{\{B_{r,t,\boldsymbol{\alpha},n}\}}\nonumber \\
&&\le \sum_{r=0}^{2^{t}-1}\sum_{\boldsymbol{\alpha} \in \{0,1\}^{n}}E\bigg|P\left (\bigcap_{i=1}^{k}\left \{\widetilde{M}_{n}(X,A_{i},\boldsymbol{\alpha})\le u_{n}(x_{i}),\widehat{M}_{n}(X,A_{i},\boldsymbol{\alpha} )\le u_{n}(y_{i})\right \}\right)\nonumber \\
&&\ \ \ \ \ - \int_{-\infty}^{+\infty}\prod_{i=1}^{k}P\left \{\widetilde{M}_{n}(Z,A_{i},\boldsymbol{\alpha})\le u_{n}(x_{i},z),\widehat{M}_{n}(Z,A_{i},\boldsymbol{\alpha})\le u_{n}(y_{i},z)\right \}\mathrm{d}\Phi(z) \bigg|I_{\{B_{r,t,\boldsymbol{\alpha},n}\}}\nonumber \\
&&+ \sum_{r=0}^{2^{t}-1}\sum_{\boldsymbol{\alpha} \in \{0,1\}^{n}}E\bigg|\int_{-\infty}^{+\infty}\prod_{i=1}^{k}\exp\left [-(d_{i}-c_{i})(\lambda g(x_{i},z)+(1-\lambda)g(y_{i},z)) \right ]\mathrm{d}\Phi(z)\nonumber \\
&&\ \ \ \ -\int_{-\infty}^{+\infty}\prod_{i=1}^{k}P\left(\widetilde{M}_{n}(Z,A_{i},\boldsymbol{\alpha})\le u_{n}(x_{i},z),\widehat{M}_{n}(Z,A_{i},\boldsymbol{\alpha})\le u_{n}(y_{i},z)\right)\mathrm{d}\Phi(z) \bigg|I_{\{B_{r,t,\boldsymbol{\alpha},n}\}}\nonumber \\
&&=\Sigma_{1}+\Sigma_{2}.
\end{eqnarray}
The first term  $\Sigma_{1}$ tends to $0$ by lemma 3.1.
For the second term, it follows from the inequality
\begin{eqnarray}
\label{15}
\bigg|\prod_{s=1}^{k}a_{s}-\prod_{s=1}^{k}b_{s}\bigg|\le \sum_{s=1}^{k}\bigg|a_{s}-b_{s} \bigg|
\end{eqnarray}
for all $a_{s}$, $b_{s}\in [0,1]$, that
\begin{eqnarray}
\label{16}
\Sigma_{2}&\le& \sum_{r=0}^{2^{t}-1}\sum_{\boldsymbol{\alpha} \in \{0,1\}^{n}}E\int_{-\infty}^{+\infty}\sum_{i=1}^{k} \bigg|\exp\left [-(d_{i}-c_{i})(\lambda g(x_{i},z)+(1-\lambda)g(y_{i},z)) \right ]\nonumber \\
&&\ \ \ \ \ \ -P\left(\widetilde{M}_{n}(Z,A_{i},\boldsymbol{\alpha})\le u_{n}(x_{i},z),\widehat{M}_{n}(Z,A_{i},\boldsymbol{\alpha})\le u_{n}(y_{i},z)\right) \bigg|\mathrm{d}\Phi(z)I_{\{B_{r,t,\boldsymbol{\alpha},n}\}}\nonumber \\
&\le& \sum_{r=0}^{2^{t}-1}\sum_{\boldsymbol{\alpha} \in \{0,1\}^{n}}E\int_{-\infty}^{+\infty}\sum_{i=1}^{k} \bigg|\exp\left [-(d_{i}-c_{i})(\lambda g(x_{i},z)+(1-\lambda)g(y_{i},z)) \right ]\nonumber \\
&&\ \ \ \ \ \ \ -\left (1-\frac{(d_{i}-c_{i})(\lambda g(x_{i},z)+(1-\lambda)g(y_{i},z))}{t}\right)^{t}\bigg|\mathrm{d}\Phi(z)I_{\{B_{r,t,\boldsymbol{\alpha},n}\}}\nonumber \\
&+&\sum_{r=0}^{2^{t}-1}\sum_{\boldsymbol{\alpha} \in \{0,1\}^{n}}E\int_{-\infty}^{+\infty}\sum_{i=1}^{k} \bigg|P\left (\widetilde{M}_{n}(Z,A_{i},\boldsymbol{\alpha})\le u_{n}(x_{i},z),\widehat{M}_{n}(Z,A_{i},\boldsymbol{\alpha})\le u_{n}(y_{i},z)\right)\nonumber \\
&&\ \ \ \ \ \ \ \ -\left (1-\frac{(d_{i}-c_{i})(\lambda g(x_{i},z)+(1-\lambda)g(y_{i},z))}{t}\right)^{t}\bigg|\mathrm{d}\Phi(z)I_{\{B_{r,t,\boldsymbol{\alpha},n}\}}\nonumber \\
&=&\Sigma_{21}+\Sigma_{22}.
\end{eqnarray}
It is easy to see that $\Sigma_{21}$ tends to $0$, as $t\rightarrow\infty$. Hence, we only need to estimate the second term $\Sigma_{22}$. We have
\begin{eqnarray}
\label{17}
\Sigma_{22}&\leq&\sum_{r=0}^{2^{t}-1}\sum_{\boldsymbol{\alpha} \in \{0,1\}^{n}}E\int_{-\infty}^{+\infty}\sum_{i=1}^{k} \bigg|P\left (\widetilde{M}_{n}(Z,A_{i},\boldsymbol{\alpha})\le u_{n}(x_{i},z),\widehat{M}_{n}(Z,A_{i},\boldsymbol{\alpha})\le u_{n}(y_{i},z)\right)\nonumber \\
&& -\prod_{s=1}^{t}P\left(\widetilde{M}_{n}(Z,K_{s,i},\boldsymbol{\alpha})\le u_{n}(x_{i},z),\widehat{M}_{n}(Z,K_{s,i},\boldsymbol{\alpha})\le u_{n}(y_{i},z)\right )\bigg|\mathrm{d}\Phi(z)I_{\{B_{r,t,\boldsymbol{\alpha},n}\}}\nonumber \\
&+&\sum_{r=0}^{2^{t}-1}\sum_{\boldsymbol{\alpha} \in \{0,1\}^{n}}E\int_{-\infty}^{+\infty}\sum_{i=1}^{k} \bigg|\prod_{s=1}^{t}P\left (\widetilde{M}_{n}(Z,K_{s,i},\boldsymbol{\alpha})\le u_{n}(x_{i},z),\widehat{M}_{n}(Z,K_{s,i},\boldsymbol{\alpha})\le u_{n}(y_{i},z)\right)\nonumber \\
&&-\prod_{s=1}^{t}\left[1-\left (\frac{r}{2^{t}}m_{i}(1-\Phi (u_{n}(x_{i},z)))+(1-\frac{r}{2^{t}})m_{i}(1-\Phi (u_{n}(y_{i},z)))  \right )  \right] \bigg|\mathrm{d}\Phi(z)I_{\{B_{r,t,\boldsymbol{\alpha},n}\}}\nonumber \\
&+&\sum_{r=0}^{2^{t}-1}\sum_{\boldsymbol{\alpha} \in \{0,1\}^{n}}E\int_{-\infty}^{+\infty}\sum_{i=1}^{k} \bigg|\prod_{s=1}^{t}\left[1-\left (\frac{r}{2^{t}}m_{i}(1-\Phi (u_{n}(x_{i},z)))+(1-\frac{r}{2^{t}})m_{i}(1-\Phi (u_{n}(y_{i},z)))\right ) \right]\nonumber \\
 &&-\left (1-\frac{(d_{i}-c_{i})(\lambda g(x_{i},z)+(1-\lambda)g(y_{i},z))}{t}\right)^{t}\bigg|\mathrm{d}\Phi(z)I_{\{B_{r,t,\boldsymbol{\alpha},n}\}}\nonumber \\
&=& \Sigma_{221}+\Sigma_{222}+\Sigma_{223}.
\end{eqnarray}
By the independence of $\{Z_{i}\}_{i\geq 1}$, it is easy to check that
\begin{eqnarray*}
\label{21}
\Sigma_{221}\leq \int_{-\infty}^{+\infty}\sum_{i=1}^{k}2m_{i}(1-\Phi(\min\{u_{n}(x_{i},z),u_{n}(y_{i},z)\}))d\Phi(z).
\end{eqnarray*}
According to the proof of Theorem 6.5.1 in Leadbetter et al. (1983), we have
$$u_{n}(x,z)=u_{n}(x+\gamma -\sqrt{2\gamma}z)+o(a_{n}^{-1}),$$ and thus,
$$\lim_{n\to \infty}(n(1-\Phi (u_{n}(x,z))))=g(x,z). $$
Consequently, we have
\begin{eqnarray}
\label{21}
\limsup_{n\rightarrow\infty}\Sigma_{221}\leq \int_{-\infty}^{+\infty}\frac{\sum_{i=1}^{k}2(d_{i}-c_{i})g(\min\{x_{i},y_{i}\},z)}{t}d\Phi(z).
\end{eqnarray}
To bound the term $\Sigma_{222}$, note that for any $0\le r\le 2^{t}-1 $, we have
\begin{eqnarray}
\label{23}
&&\sum_{j\in K_{s,i}}\alpha _{j}(\Phi (u_{n}(x_{i},z))-\Phi (u_{n}(y_{i},z)))+[1-m_{i}(1-\Phi (u_{n}(y_{i},z)))]\nonumber \\
&&=\left [ 1-\frac{m_{i}r}{2^{t}} (1-\Phi (u_{n}(x_{i},z))) -m_{i}(1-\frac{r}{2^{t}})(1-\Phi (u_{n}(y_{i},z)))\right ]\nonumber \\
&&\ \ \ +\left [ \frac{\sum_{j\in K_{s,i}}\alpha _{j}}{m_{i}}-\frac{r}{2^{t}}\right ]m_{i}(\Phi (u_{n}(x_{i},z))-\Phi (u_{n}(y_{i},z)))\nonumber \\
&&\le P\left \{\widetilde{M}_{n}(Z,K_{s,i},\boldsymbol{\alpha})\le u_{n}(x_{i},z),\widehat{M}_{n}(Z,K_{s,i},\boldsymbol{\alpha})\le u_{n}(y_{i},z)\right \}\nonumber \\
&&\le m_{i}\sum_{j=2}^{m_{i}}P[A_{s1},A_{sj}]+\left [ 1-\frac{m_{i}r}{2^{t}} (1-\Phi (u_{n}(x_{i},z))) -m_{i}(1-\frac{r}{2^{t}})(1-\Phi (u_{n}(y_{i},z)))\right ]\nonumber \\
&&\ \ \ +\left [ \frac{\sum_{j\in K_{s,i}}\alpha _{j}}{m_{i}}-\frac{r}{2^{t}}\right ]m_{i}(\Phi (u_{n}(x_{i},z))-\Phi (u_{n}(y_{i},z))),
\end{eqnarray}
where $A_{sj}=\{Z_{(s-1)m_{i}+j}>\min\{u_{n}(x_{i},z),u_{n}(y_{i},z), i=1,2,\ldots,k\}\} $, $j\in \{1,2,\cdots ,m_{i}\}$.
Thus, combining Lemma 3.2 with \eqref{15} and \eqref{23}, we get
\begin{eqnarray}
\label{24}
\Sigma_{222}&\le& \sum_{r=0}^{2^{t}-1}\sum_{\boldsymbol{\alpha} \in \{0,1\}^{n}}E\int_{-\infty}^{+\infty}\sum_{i=1}^{k}\sum_{s=1}^{t}\bigg|\left [ \frac{\sum_{j\in K_{s,i}}\alpha _{j}}{m_{i}}-\frac{r}{2^{t}}\right ]m_{i}(\Phi (u_{n}(x_{i},z))-\Phi (u_{n}(y_{i},z)))\nonumber \\
&&\ \ \ \  +m_{i}\sum_{j=2}^{m_{i}}P[A_{s1},A_{sj}]  \bigg|\mathrm{d}\Phi(z)I_{\{B_{r,t,\boldsymbol{\alpha},n}\}}\nonumber \\
&\le& \sum_{r=0}^{2^{t}-1}E\int_{-\infty}^{+\infty}\sum_{i=1}^{k}\sum_{s=1}^{t}m_{i}\bigg| \frac{\sum_{j\in K_{s,i}}\varepsilon _{j}}{m_{i}}-\frac{r}{2^{t}}\bigg|\bigg|(\Phi (u_{n}(x_{i},z))-\Phi (u_{n}(y_{i},z)))\bigg|\mathrm{d}\Phi(z)I_{\{B_{r,t}\}}\nonumber \\
&&\ \ \ \ +\int_{-\infty}^{+\infty}\sum_{i=1}^{k}\sum_{s=1}^{t}[m_{i}^{2}(1-\Phi \{\min(u_{n}(x_{i},z),u_{n}(y_{i},z))\})^{2}] \mathrm{d}\Phi(z)\nonumber \\
&\le&\sum_{i=1}^{k}\sum_{s=1}^{t}\left(E\bigg|\frac{\sum_{j\in K_{s,i}}\varepsilon _{j}}{m_{i}}-\lambda \bigg|+\sum_{r=0}^{2^{t}-1}E\bigg|\lambda-\frac{r}{2^{t}}\bigg|I_{\{B_{r,t}\}}\right)\nonumber \\
&&\ \ \ \ \times\int_{-\infty}^{+\infty}m_{i}\bigg|(\Phi (u_{n}(x_{i},z))-\Phi (u_{n}(y_{i},z)))\bigg|\mathrm{d}\Phi(z) \nonumber \\
&&\ \ \ \ +\int_{-\infty}^{+\infty}\sum_{i=1}^{k}\sum_{s=1}^{t}[m_{i}^{2}(1-\Phi \{\min(u_{n}(x_{i},z),u_{n}(y_{i},z))\})^{2}] \mathrm{d}\Phi(z)\nonumber \\
&\le&\sum_{i=1}^{k}\sum_{s=1}^{t}\left \{ \left [2(2s-1)\left ( d\left ( \frac{S_{m_{i}s+\lfloor nc_{i}\rfloor}}{m_{i}s+\lfloor nc_{i}\rfloor} ,\lambda \right )+d\left ( \frac{S_{m_{i}(s-1)+\lfloor nc_{i}\rfloor}}{m_{i}(s-1)+\lfloor nc_{i}\rfloor} ,\lambda \right ) \right )  +\frac{1}{2^{t}}\right ]\right. \nonumber \\
&&\left.\ \ \ \ \times \int_{-\infty}^{+\infty}m_{i}\bigg|(\Phi (u_{n}(x_{i},z))-\Phi (u_{n}(y_{i},z)))\bigg| \mathrm{d}\Phi(z) \right \}\nonumber \\
&&\ \ \ \ +\int_{-\infty}^{+\infty}\sum_{i=1}^{k}\sum_{s=1}^{t}[m_{i}^{2}(1-\Phi \{\min(u_{n}(x_{i},z),u_{n}(y_{i},z))\})^{2}] \mathrm{d}\Phi(z).
\end{eqnarray}
Since $\lim_{m_{i}\to \infty} d\left ( \frac{S_{m_{i}s+\lfloor nc_{i}\rfloor}}{m_{i}s+\lfloor nc_{i}\rfloor} ,\lambda \right )=0 $, we obtain
\begin{eqnarray}
\label{25}
\limsup_{n \to \infty} \Sigma_{222}&\le& \sum_{i=1}^{k}\frac{1}{2^{t}}\int_{-\infty}^{+\infty}(d_{i}-c_{i})|g(y_{i},z)-g(x_{i},z)|\mathrm{d}\Phi(z)\nonumber \\
&&+\int_{-\infty}^{+\infty}\frac{\sum_{i=1}^{k}(d_{i}-c_{i})^{2}(g(\min\{x_{i},y_{i}\},z)) ^{2}}{t}\mathrm{d}\Phi(z).
\end{eqnarray}
For $\Sigma_{223}$, using \eqref{15} again, we get
\begin{eqnarray*}
\Sigma_{223}&\le& \sum_{r=0}^{2^{t}-1}E\int_{-\infty}^{+\infty}\sum_{i=1}^{k}\sum_{s=1}^{t}\bigg|\frac{(d_{i}-c_{i})(\lambda g(x_{i},z)+(1-\lambda)g(y_{i},z))}{t}\nonumber \\
&& -\frac{\frac{r}{2^{t}}n(d_{i}-c_{i})(1-\Phi (u_{n}(x_{i},z)))+(1-\frac{r}{2^{t}})n(d_{i}-c_{i})(1-\Phi (u_{n}(y_{i},z)))}{t} \bigg|\mathrm{d}\Phi(z)I_{\{B_{r,t}\}}.
\end{eqnarray*}
Thus,
\begin{eqnarray}
\label{26}
\limsup_{n\to \infty} \Sigma_{223}&\le& \int_{-\infty}^{+\infty}\sum_{i=1}^{k}\sum_{s=1}^{t}\sum_{r=0}^{2^{t}-1}E\left (\frac{\bigg|\lambda -\frac{r}{2^{t}}\bigg|(d_{i}-c_{i})(g(x_{i},z)+g(y_{i},z))}{t}\right )\mathrm{d}\Phi(z)I_{\{B_{r,t}\}}\nonumber \\
&\le& \sum_{i=1}^{k} \int_{-\infty }^{+\infty}\frac{(d_{i}-c_{i})(g(x_{i},z)+g(y_{i},z))}{2^{t}} \mathrm{d}\Phi(z).
\end{eqnarray}
From \eqref{17}, \eqref{21}, \eqref{25} and \eqref{26} we have
\begin{eqnarray}
\label{27}
\limsup_{n\to \infty} \Sigma_{22}&\le&\int_{-\infty}^{+\infty}\frac{\sum_{i=1}^{k}2(d_{i}-c_{i})g(\min\{x_{i},y_{i}\},z)}{t}d\Phi(z)\nonumber \\
&&+\sum_{i=1}^{k}\frac{1}{2^{t}}\int_{-\infty}^{+\infty}(d_{i}-c_{i})|g(y_{i},z)-g(x_{i},z)|\mathrm{d}\Phi(z)\nonumber \\
&&+\int_{-\infty}^{+\infty}\frac{\sum_{i=1}^{k}(d_{i}-c_{i})^{2}(g(\min\{x_{i},y_{i}\},z)) ^{2}}{t}\mathrm{d}\Phi(z)\nonumber \\
&&+\sum_{i=1}^{k} \int_{-\infty }^{+\infty}\frac{(d_{i}-c_{i})(g(x_{i},z)+g(y_{i},z))}{2^{t}} \mathrm{d}\Phi(z).
\end{eqnarray}
Now if we take $t\rightarrow\infty$ of both sides of (\ref{27}), we obtain $\Sigma_{22}\rightarrow0$, which completes the proof.

The following lemma is taken from Wi\'{s}niewski (1994).

\textbf{Lemma 3.4}. {\sl Suppose that $(M_{n1},M_{n2},\cdots, M_{nd} )$ is a vector point process of exceedances of the level $w_{nk} $ by random vector $\{\eta_{n1},\eta_{n2},\ldots,\eta_{nd}\}_{n\ge 1}$,  i.e.,
$$M_{nk}(B)=\sum_{\frac{j}{n}\in B}I_{\{\eta_{nk}>w_{nk}\}}, k=1,2,\cdots, d, $$
for any Borel set $B\subset (0,1]$. Suppose that

(i) $M_{nk}$ converges in distribution to the point process $M_{k}$ and $M_{k}$ is simple point process.

(ii) For any Borel sets $B_{k}\subset(0,1], k=1,2,\cdots, d$, $\lim_{n \to \infty}P(M_{nk}(B_{k})=0,\ \ k=1,2\cdots, d )$ exists. \\
Then the vector point process $(M_{n1},M_{n2},\cdots ,M_{nd} )$ converges in distribution to $(M_{1},M_{2},\cdots,M_{d} )$  on $(0,1]^{d}$ as $n \rightarrow \infty$.}

\section{Proof of main results}
In this section, we give the proofs of our main results.

\textbf{The proof of Theorem 2.1}. Firstly,  we prove that $\widetilde{N}_{n}^{x} $ and $\widehat{N}_{n}^{y}$ converge in distribution to $\widetilde{N}^{x}$ and $\widehat{N}^{y}$, respectively. From Kallenberg's Theorem (see e.g., Kallenberg (1975), Leadbetter et al. (1983)), it suffices to show that as $n\rightarrow\infty $,

(a) $E[\widetilde{N}_{n}^{x}(B_{1})]\rightarrow E[\widetilde{N}^{x}(B_{1})]$, $E[\widehat{N}_{n}^{y}(B_{1})]\rightarrow E[\widehat{N}^{y}(B_{1})]$, where $B_{1}=(c,d]$, $0< c<d\leq1$;

(b) $P(\widetilde{N}_{n}^{x}(B_{2})=0)\rightarrow P(\widetilde{N}^{x}(B_{2})=0)$, $P(\widehat{N}_{n}^{y}(B_{2})=0)\rightarrow P(\widehat{N}^{y}(B_{2})=0)$, where $B_{2}=\bigcup_{i=1}^{k}(c_{i},d_{i}] $, $0< c_{1}<d_{1}\le c_{2}<\cdots <c_{k}<d_{k}\le 1$.
\\To show assertion (a),  let $A=\left \{\lfloor nc\rfloor+1,\lfloor nc\rfloor+2,\cdots ,\lfloor nd\rfloor\right \}$. According to the conditions of the Theorem 2.1, we have
$$\frac{S_{n}(A)}{\lfloor nd\rfloor-\lfloor nc\rfloor}\overset{p}{\longrightarrow} \lambda, \ \ \ as \ \ n\to \infty, $$
which combined with the dominated convergence theorem yields
$$ E\left (\frac{S_{n}(A)}{\lfloor nd\rfloor-\lfloor nc\rfloor}\right ) \rightarrow E(\lambda), \ \ \ as \ \ n\to \infty,   $$ i.e.,
\begin{eqnarray}
\label{28}
\ \ \ \frac{\sum_{j\in A}P(\varepsilon _{j}=1 )}{n(d-c)}  \rightarrow E(\lambda), \ \ \ as \ \ n\to \infty.
\end{eqnarray}
Thus, we also have
\begin{eqnarray}
\label{29}
\ \ \ \frac{\sum_{j\in A}P(\varepsilon _{j}=0 )}{n(d-c)}  \rightarrow E(1-\lambda), \ \ \ as \ \ n\to \infty.
\end{eqnarray}
Using the facts (\ref{28}) and $\lim_{n \to \infty} n(1-\Phi (u_{n}(x))) =\exp \left (-x \right )$, we have
\begin{eqnarray}
\label{30}
E[\widetilde{N}_{n}^{x}(B_{1})]&=&E\left[\sum_{\frac{j}{n}\in (c,d]}I_{\{X_{j}>u_{n}(x),\varepsilon _{j}=1\}}\right ]\nonumber \\
&=&(1-\Phi(u_{n}(x)))\sum_{j\in A}P(\varepsilon _{j}=1)\nonumber \\
&\rightarrow& (d-c)E(\lambda) e^{-x}.
\end{eqnarray}
Noting that $\widetilde{N}^{x}$ is a conditional Poisson process with conditional intensity $\lambda g(x,\xi)$, we get
\begin{eqnarray}
\label{31}
E[\widetilde{N}^{x}(B_{1})]&=&E\{E[\widetilde{N}^{x}(B_{1})|(\lambda ,\xi)]\}\nonumber \\
&=&E((d-c)\lambda g(x,\xi))\nonumber \\
&=&(d-c)E(\lambda)E(g(x,\xi))\nonumber \\
&=&(d-c)E(\lambda) e^{-x}.
\end{eqnarray}
Similarly,
$$E[\widehat{N}_{n}^{y}(B_{1})]\rightarrow(d-c)E((1-\lambda)) e^{-y}=E[\widehat{N}^{y}(B_{1})].$$
Thus, the first assertion holds.\\
To prove the second assertion (b), using Lemma 3.3 for the case $x_{i}\equiv x$ and $y_{i}\equiv y$, $i=1,2,\ldots,k$, we have
\begin{eqnarray}
\label{32}
&&\lim_{n \to \infty} P(\widetilde{N}_{n}^{x}(B_{2})=0,\widehat{N}_{n}^{y}(B_{2})=0)\nonumber \\
&&\ \ \ =\lim_{n \to \infty}P\left (\bigcap_{i=1}^{k}\left \{\widetilde{M}_{n}(X,A_{i},\boldsymbol
{\varepsilon})\le u_{n}(x),\widehat{M}_{n}(X,A_{i},\boldsymbol
{\varepsilon})\le u_{n}(y)\right \}\right)\nonumber \\
&&\ \ \ =E\int_{-\infty}^{+\infty}\prod_{i=1}^{k}\exp\left [-(d_{i}-c_{i})(\lambda g(x,z)+(1-\lambda)g(y,z)) \right ]\mathrm{d}\Phi(z).
\end{eqnarray}
On the one hand, letting $y\to \infty$ and $x\to \infty$, we obtain
\begin{eqnarray}
\label{33}
\lim_{n \to \infty} P(\widetilde{N}_{n}^{x}(B_{2})=0)=E\int_{-\infty}^{+\infty}\prod_{i=1}^{k}\exp \left [-(d_{i}-c_{i})(\lambda g(x,z))\right ]\mathrm{d}\Phi(z)
\end{eqnarray}
and
\begin{eqnarray}
\label{34}
\lim_{n\to \infty} P(\widehat{N}_{n}^{y}(B_{2})=0)=E\int_{-\infty}^{+\infty}\prod_{i=1}^{k}\exp \left [-(d_{i}-c_{i})((1-\lambda) g(y,z))\right ]\mathrm{d}\Phi(z),
\end{eqnarray}
respectively. On the other hand,  since $\widetilde{N}^{x}$ is a conditional Poisson process with conditional intensity $\lambda g(x,\xi)$, we have
\begin{eqnarray}
\label{35}
P(\widetilde{N}^{x}(B_{2})=0)&=&P\left ( \bigcap_{i=1}^{k}\{\widetilde{N}^{x}((c_{i},d_{i}])=0\}\right )\nonumber \\
&=&\int_{0}^{1}\int_{-\infty}^{+\infty}\prod_{i=1}^{k} P(\widetilde{N}^{x}(c_{i},d_{i}]=0|\lambda =e,\xi =z)\mathrm{d}F_{\lambda}(e)\mathrm{d}\Phi(z)\nonumber \\
&=&E\int_{-\infty}^{+\infty}\prod_{i=1}^{k}\exp \left [-(d_{i}-c_{i})(\lambda g(x,z))\right ]\mathrm{d}\Phi(z),
\end{eqnarray}
where $F_{\lambda}(e)=P(\lambda \le e) $.
Hence, $P(\widetilde{N}_{n}^{x}(B_{2})=0)\rightarrow P(\widetilde{N}^{x}(B_{2})=0)$, as $n\rightarrow\infty$.
Similarly, we have
$$P(\widehat{N}_{n}^{y}(B_{2})=0)\rightarrow P(\widehat{N}^{y}(B_{2})=0),$$
 as $n\rightarrow\infty$, which proves (b).
\\Secondly, we show that for any Borel sets $C_{i}\subset (0,1]$, $i=1,\ldots,k+l$,
\begin{eqnarray}
\label{36-0}
&&\lim_{n \to \infty} P\left ( \widetilde{N}_{n}^{x_{1}}(C_{1})=0,\ldots,\widetilde{N}_{n}^{x_{k}}(C_{k})=0,\widehat{N}_{n}^{y_{1}}(C_{k+1})=0, \ldots ,\widehat{N}_{n}^{y_{l}}(C_{k+l})=0\right )
\end{eqnarray}
exists. For the simplicity, we only deal with the case $k=2$ and $l=2$, since the general case can be done similarly.
We first deal with the case $x_{1}\leq x_{2}$ and $y_{1}\leq y_{2}$.
Put $nC_{i} =\{nx,x\in C_{i}\}$ and let $C_{5}=C_{2}\setminus C_{1}$, $C_{6}=C_{4}\setminus C_{3}$. We have
\begin{eqnarray*}
&&\lim_{n \to \infty} P\left ( \widetilde{N}_{n}^{x_{1}}(C_{1})=0, \widetilde{N}_{n}^{x_{2}}(C_{2})=0, \widehat{N}_{n}^{y_{1}}(C_{3})=0 , ,\widehat{N}_{n}^{y_{2}}(C_{4})=0\right ) \\
&&\ \ \ \ \  =\lim_{n \to \infty}P\left ( \widetilde{M}_{C_{1}}\le u_{n}(x_{1}),\widetilde{M}_{C_{2}}\le u_{n}(x_{2}),\widehat{M}_{C_{3}}\le u_{n}(y_{1}),\widehat{M}_{C_{4}}\le u_{n}(y_{2})\right) \\
&&\ \ \ \ \  =\lim_{n \to \infty}P\left ( \widetilde{M}_{C_{1}}\le u_{n}(x_{1}),\widetilde{M}_{C_{5}}\le u_{n}(x_{2}),\widehat{M}_{C_{3}}\le u_{n}(y_{1}),\widehat{M}_{C_{6}}\le u_{n}(y_{2})\right),
\end{eqnarray*}
where $\widetilde{M}_{C_{i}}=\widetilde{M}_{n}(X,nC_{i},\boldsymbol{\varepsilon})$ and $\widehat{M}_{C_{i}}=\widehat{M}_{n}(X,nC_{i},\boldsymbol{\varepsilon})$.
Note that $C_{1}\bigcap C_{5}=\varnothing$ and $C_{3}\bigcap C_{6}=\varnothing$.
It follows now from Lemma 3.3 that
\begin{eqnarray}
\label{36}
&&\lim_{n \to \infty} P\left ( \widetilde{N}_{n}^{x_{1}}(C_{1})=0, \widetilde{N}_{n}^{x_{2}}(C_{2})=0, \widehat{N}_{n}^{y_{1}}(C_{3})=0 , ,\widehat{N}_{n}^{y_{2}}(C_{4})=0\right ) \nonumber \\
&&=\lim_{n \to \infty}P\left ( \widetilde{M}_{C_{1}}\le u_{n}(x_{1}),\widetilde{M}_{C_{5}}\le u_{n}(x_{2}),\widehat{M}_{C_{3}}\le u_{n}(y_{1}),\widehat{M}_{C_{6}}\le u_{n}(y_{2})\right)\nonumber \\
&&=E\int_{-\infty }^{+\infty}\exp\left (-m(C_{1})\lambda g(x_{1},z)-m(C_{3})(1-\lambda)g(y_{1},z)\right) \nonumber \\
&&\ \ \ \ \ \ \times\exp\left( -m(C_{2}\setminus C_{1})\lambda g(x_{2},z)-m(C_{4}\setminus C_{3})(1-\lambda) g(y_{2},z)\right)\mathrm{d}\Phi(z).
\end{eqnarray}
By the same arguments as the above, we have, for $x_{1}\leq x_{2}$ and $y_{1}> y_{2}$,
\begin{eqnarray}
\label{37}
&&\lim_{n \to \infty} P\left ( \widetilde{N}_{n}^{x_{1}}(C_{1})=0, \widetilde{N}_{n}^{x_{2}}(C_{2})=0, \widehat{N}_{n}^{y_{1}}(C_{3})=0 , ,\widehat{N}_{n}^{y_{2}}(C_{4})=0\right ) \nonumber \\
&&=E\int_{-\infty }^{+\infty}\exp\left (-m(C_{1})\lambda g(x_{1},z)-m(C_{4})(1-\lambda)g(y_{2},z)\right) \nonumber \\
&&\ \ \ \ \ \ \times\exp\left( -m(C_{2}\setminus C_{1})\lambda g(x_{2},z)-m(C_{3}\setminus C_{4})(1-\lambda) g(y_{1},z)\right)\mathrm{d}\Phi(z);
\end{eqnarray}
for $x_{1}> x_{2}$ and $y_{1}\leq y_{2}$,
\begin{eqnarray}
\label{37-1}
&&\lim_{n \to \infty} P\left ( \widetilde{N}_{n}^{x_{1}}(C_{1})=0, \widetilde{N}_{n}^{x_{2}}(C_{2})=0, \widehat{N}_{n}^{y_{1}}(C_{3})=0 , ,\widehat{N}_{n}^{y_{2}}(C_{4})=0\right ) \nonumber \\
&&=E\int_{-\infty }^{+\infty}\exp\left (-m(C_{2})\lambda g(x_{2},z)-m(C_{3})(1-\lambda)g(y_{1},z)\right) \nonumber \\
&&\ \ \ \ \ \ \times\exp\left( -m(C_{1}\setminus C_{2})\lambda g(x_{1},z)-m(C_{4}\setminus C_{3})(1-\lambda) g(y_{2},z)\right)\mathrm{d}\Phi(z)
\end{eqnarray}
and  for $x_{1}>x_{2}$ and $y_{1}> y_{2}$,
\begin{eqnarray}
\label{37-2}
&&\lim_{n \to \infty} P\left ( \widetilde{N}_{n}^{x_{1}}(C_{1})=0, \widetilde{N}_{n}^{x_{2}}(C_{2})=0, \widehat{N}_{n}^{y_{1}}(C_{3})=0 , ,\widehat{N}_{n}^{y_{2}}(C_{4})=0\right ) \nonumber \\
&&=E\int_{-\infty }^{+\infty}\exp\left (-m(C_{2})\lambda g(x_{2},z)-m(C_{4})(1-\lambda)g(y_{2},z)\right) \nonumber \\
&&\ \ \ \ \ \ \times\exp\left( -m(C_{1}\setminus C_{2})\lambda g(x_{1},z)-m(C_{3}\setminus C_{4})(1-\lambda) g(y_{1},z)\right)\mathrm{d}\Phi(z).
\end{eqnarray}
Therefore, by Lemma 3.4, we prove that $\left ( \widetilde{N}_{n}^{x_{1}},\ldots,\widetilde{N}_{n}^{x_{k}}, \widehat{N}_{n}^{y_{1}},\ldots,\widehat{N}_{n}^{y_{l}}\right ) $ converges in distribution to the point process $\left ( \widetilde{N}^{x_{1}},\ldots, \widetilde{N}^{x_{k}}, \widehat{N}^{y_{1}},\ldots ,\widehat{N}^{y_{l}}\right ) $ on $\left ( 0,1 \right ] ^{k+l} $.

\textbf{Proof of Corollary 2.1}.   By Theorem 2.1, we have
\begin{eqnarray}
\label{38}
&&P\left ( \widetilde{N}^{x}(B)=k_{1} ,\widehat{N}^{x}(B)=k_{2} , \widetilde{N}^{y}(B)=k_{3} ,\widehat{N}^{y}(B)=k_{4}\right )\nonumber \\
&&=\int_{-\infty}^{+\infty}\int_{0}^{1}P\left (\widetilde{N}^{x}(B)=k_{1}, \widetilde{N}^{y}(B)=k_{3}|\lambda =e,\xi =z \right )\nonumber \\
&&\ \ \ \ \ \ \times P\left ( \widehat{N}^{x}(B)=k_{2},\widehat{N}^{y}(B)=k_{4}|\lambda =e,\xi =z \right )\mathrm{d}\Phi(z)\mathrm{d}F_{\lambda}(e).
\end{eqnarray}
For $x>y$, applying Theorem 5.6.1 in Leadbetter et al. (1983), we have
\begin{eqnarray}
\label{39}
&&P\left (\widetilde{N}^{x}(B)=k_{1}, \widetilde{N}^{y}(B)=k_{3}|\lambda =e,\xi =z \right )\nonumber \\
&&\ \ =\frac{[em(B)(g(y,z)-g(x,z))]^{k_{3}-k_{1}}}{(k_{3}-k_{1})!}\frac{[em(B)g(x,z)]^{k_{1}} }{k_{1}!} \exp(-em(B)g(y,z)).
\end{eqnarray}
Similarly,
\begin{eqnarray}
\label{40}
&&P\left (\widehat{N}^{x}(B)=k_{2}, \widehat{N}^{y}(B)=k_{4}|\lambda =e,\xi =z \right )\nonumber \\
&&\ \ =\frac{[(1-e)m(B)(g(y,z)-g(x,z))]^{k_{2}-k_{4}}}{(k_{2}-k_{4})!}\frac{[(1-e)m(B)g(x,z)]^{k_{2}} }{k_{2}!} \exp(-(1-e)m(B)g(y,z)).\nonumber \\
\end{eqnarray}
Hence, plugging (\ref{39}) and(\ref{40}) into (\ref{38}), we obtain the desired result (\ref{2.3}). By the same arguments, we can prove (\ref{2.4}) also holds.

\textbf{Proof of Theorem  2.2}. (i). First, note that
$$\{\widetilde{M}_{n}^{(k)}\le u_{n}(x)\}=\{\widetilde{N}_{n}^{x}((0,1]) \le k-1 \}$$
and
$$\{\widehat{M}_{n}^{(l)}\le u_{n}(y)\}=\{\widehat{N}_{n}^{y}((0,1]) \le l-1 \}.$$
Thus, by Theorem 2.1, we have
\begin{eqnarray}
\label{41}
&&\lim_{n \to \infty}P\left \{ \widetilde{M}_{n}^{(k)}\le u_{n}\left ( x \right ),\widehat{M}_{n} ^{(l)} \le u_{n}\left ( y \right ) \right \}\nonumber \\
&&=\lim_{n \to \infty}P\left (\widetilde{N}_{n}^{x}((0,1]) \le k-1,\widehat{N}_{n}^{y}((0,1]) \le l-1 \right )\nonumber \\
&&=P\left (\widetilde{N}^{x}((0,1]) \le k-1,\widehat{N}^{y}((0,1]) \le l-1 \right )\nonumber \\
&&=\int_{0}^{1}\int_{-\infty}^{+\infty}P\left (\widetilde{N}^{x}((0,1])\le k-1|\lambda =e,\xi =z\right )\nonumber \\
&& \ \ \ \ \ \ \ \ \ \times P\left (\widehat{N}^{y}((0,1])\le l-1|\lambda =e,\xi =z\right )\mathrm{d}F_{\lambda}(e)\mathrm{d}\Phi (z)\nonumber \\
&&=\sum_{s=0}^{k-1}\sum_{t=0}^{l-1}E\int_{-\infty}^{+\infty}\left \{\frac{[\lambda g(x,z)]^{s} }{s!}\exp(-\lambda g(x,z))\right.\nonumber \\
&&\ \ \ \ \ \ \ \ \ \  \ \left. \times\frac{[(1-\lambda)g(y,z)]^{t}}{t!}\exp(-(1-\lambda)g(y,z)) \right \}\mathrm{d}\Phi(z),
\end{eqnarray}
which proves the (i).

(ii). By the same way as (i), we have
\begin{eqnarray}
\label{42}
&&\lim_{n\to \infty}P\left ( \widetilde{M}_{n}^{(k)}\le u_{n}(x),M_{n}^{(m)}\le u_{n}(y)\right )\nonumber \\
&&\ \ \ \ \ \ =\lim_{n\to \infty}P\left ( \widetilde{N}_{n}^{x}((0,1])\le k-1, \widetilde{N}_{n}^{y}((0,1])+\widehat{N}_{n}^{y}((0,1])\le m-1 \right )\nonumber \\
&&\ \ \ \ \ \ =P\left ( \widetilde{N}^{x}((0,1])\le k-1, \widetilde{N}^{y}((0,1])+\widehat{N}^{y}((0,1])\le m-1 \right )\nonumber \\
&&\ \ \ \ \ \ =\int_{0}^{1}\int_{-\infty}^{+\infty}\sum_{t=0}^{m-1}P\left (\widehat{N}^{y}((0,1])= t|\lambda =e,\xi =z\right )\nonumber \\
&&\ \ \ \ \ \ \ \ \times P\left (\widetilde{N}^{x}((0,1])\le k-1,\widetilde{N}^{y}((0,1])\le m-1-t|\lambda =e,\xi =z\right ) \mathrm{d}F_{\lambda}(e)\mathrm{d}\Phi (z)\nonumber \\
&&\ \ \ \ \ \ =\int_{0}^{1}\int_{-\infty}^{+\infty}\sum_{t=0}^{m-1}\frac{[(1-e)g(y,z)]^{t}\exp(-(1-e)g(y,z)) }{t!}\nonumber \\
&&\ \ \ \ \ \ \ \ \times P\left (\widetilde{N}^{x}((0,1])\le k-1,\widetilde{N}^{y}((0,1])\le m-1-t|\lambda =e,\xi =z\right ) \mathrm{d}F_{\lambda}(e)\mathrm{d}\Phi (z).
\end{eqnarray}
For $x<y$, using (\ref{39}), we have
\begin{eqnarray}
\label{43}
&&P\left (\widetilde{N}^{x}((0,1])\le k-1,\widetilde{N}^{y}((0,1])\le m-1-t|\lambda =e,\xi =z\right )\nonumber \\
&&\ \ \ \ \ =\sum_{i=0}^{m-t-1}\sum_{j=i}^{k-1}P\left (\widetilde{N}^{x}((0,1])=j,\widetilde{N}^{y}((0,1])=i|\lambda =e,\xi =z \right )\nonumber \\
&&\ \ \ \ \ =\sum_{i=0}^{m-t-1}\sum_{j=i}^{k-1} \frac{[eg(y,z)]^{i}}{i!}\frac{[e(g(x,z)-g(y,z))]^{j-i}}{(j-i)!}\exp(-eg(x,z)).
\end{eqnarray}
Plugging (\ref{43}) into (\ref{42}), we get the desired result (\ref{2.6}).
For $x\geq y$, similarly, we have
\begin{eqnarray}
\label{44}
&&P\left (\widetilde{N}^{x}((0,1])\le k-1,\widetilde{N}^{y}((0,1])\le m-1-t|\lambda =e,\xi =z\right )\nonumber \\
&&\ \ \ \ \ =\sum_{i=0}^{m-t-1}\sum_{j=0}^{\min\{i,k-1\}}P\left (\widetilde{N}^{x}((0,1])=j,\widetilde{N}^{y}((0,1])=i|\lambda =e,\xi =z \right )\nonumber \\
&&\ \ \ \ \ =\sum_{i=0}^{m-t-1}\sum_{j=0}^{\min\{i,k-1\}}\frac{[eg(x,z)]^{j}}{j!}\frac{[e(g(y,z)-g(x,z))]^{i-j}}{(i-j)!}\exp(-eg(y,z)).
\end{eqnarray}
Plugging (\ref{44}) into (\ref{42}), we get the desired result (\ref{2.7}).

\textbf{Proof of Theorem 2.3}. (i). Recall that $\mathbb{N}_{n}=\{1,2,\ldots,n\}$ and define $\mathbb{N}_{ns}=\{1,2,\ldots,\lfloor ns\rfloor\}$.
Note that
\begin{eqnarray}
\label{c231}
&&P\left(\frac{1}{n}\widetilde{L}_{n}\leq s, \frac{1}{n}\widehat{L}_{n}\leq t, \widetilde{M}_{n}\leq u_{n}(x), \widehat{M}_{n}\leq u_{n}(y)\right)\nonumber\\
&&=P\left(\widetilde{M}_{n}(X,\mathbb{N}_{n}\setminus \mathbb{N}_{ns}, \boldsymbol{\varepsilon}) \leq \widetilde{M}_{n}(X,\mathbb{N}_{ns}, \boldsymbol{\varepsilon}) \leq u_{n}(x),\right.\nonumber\\
 &&\ \ \ \ \ \ \ \ \left.\widehat{M}_{n}(X,\mathbb{N}_{n}\setminus \mathbb{N}_{nt}, \boldsymbol{\varepsilon}) \leq \widehat{M}_{n}(X,\mathbb{N}_{nt}, \boldsymbol{\varepsilon}) \leq u_{n}(y)\right)=: P_{n}(s,t,x,y).
\end{eqnarray}
Define $\mathbb{R}_{st}=(s,t]$ for $0\leq s<t\leq 1$. By Theorem 2.1,  we have
\begin{eqnarray}
\label{c232}
&&P\left(\widetilde{M}_{n}(X,\mathbb{N}_{n}\setminus \mathbb{N}_{ns}, \boldsymbol{\varepsilon}) \leq u_{n}(x_{1}),  \widetilde{M}_{n}(X,\mathbb{N}_{ns}, \boldsymbol{\varepsilon}) \leq u_{n}(x_{2}),\right.\nonumber\\
 &&\ \ \ \ \ \ \ \ \left.\widehat{M}_{n}(X,\mathbb{N}_{n}\setminus \mathbb{N}_{nt}, \boldsymbol{\varepsilon}) \leq u_{n}(x_{3}), \widehat{M}_{n}(X,\mathbb{N}_{nt}, \boldsymbol{\varepsilon}) \leq u_{n}(x_{4})\right)\nonumber\\
 &&=P\left(\widetilde{N}_{n}^{x_{1}}(\mathbb{R}_{s1})=0, \widetilde{N}_{n}^{x_{2}}(\mathbb{R}_{0s})=0,
 \widehat{N}_{n}^{x_{3}}(\mathbb{R}_{t1})=0, \widehat{N}_{n}^{x_{4}}(\mathbb{R}_{0t})=0\right)\nonumber\\
 &&\rightarrow P\left(\widetilde{N}^{x_{1}}(\mathbb{R}_{s1})=0, \widetilde{N}^{x_{2}}(\mathbb{R}_{0s})=0,
 \widehat{N}^{x_{3}}(\mathbb{R}_{t1})=0, \widehat{N}^{x_{4}}(\mathbb{R}_{0t})=0\right)\nonumber\\
&&= \int_{0}^{1}\int_{-\infty}^{+\infty}P\left (\widetilde{N}^{x_{1}}(\mathbb{R}_{s1})=0, \widetilde{N}^{x_{2}}( \mathbb{R}_{0s})=0|\lambda =e,\xi =z\right )\nonumber \\
&& \ \ \ \ \ \ \ \ \ \times P\left (\widehat{N}^{x_{3}}(\mathbb{R}_{t1})=0,\widehat{N}^{x_{4}}(\mathbb{R}_{0t})=0 |\lambda =e,\xi =z\right )\mathrm{d}F_{\lambda}(e)\mathrm{d}\Phi (z)\nonumber \\
&&=\int_{0}^{1}\int_{-\infty}^{+\infty} \exp\left(-(1-s)eg(x_{1},z)-seg(x_{2},z)\right)\nonumber \\
&&\ \ \ \ \ \ \ \ \ \ \times\exp\left(-(1-t)(1-e)g(x_{3},z)-t(1-e)g(x_{4},z)\right)\mathrm{d}F_{\lambda}(e)\mathrm{d}\Phi (z)\nonumber \\
&&=: H(x_{1},x_{2},x_{3},x_{4}).
\end{eqnarray}
Let $h(x_{1},x_{2},x_{3},x_{4})=\frac{\partial^{4}H(x_{1},x_{2},x_{3},x_{4})}{\partial x_{1}\partial x_{2}\partial x_{3}\partial x_{4}}$ be the joint density function of random vector $(Y_{1}, Y_{2}, Y_{3}, Y_{4})$, then by some direct computations
\begin{eqnarray}
\label{c233}
 \lim_{n\rightarrow\infty}P_{n}(s,t,x,y)&=&P(Y_{1}\leq Y_{2}\leq x, Y_{3}\leq Y_{4}\leq y)\nonumber \\
 &=&\iiiint\limits_{x_{1}\leq x_{2}\leq x, x_{3}\leq x_{4}\leq y}h(x_{1},x_{2},x_{3},x_{4})
 \mathrm{d}x_{1}\mathrm{d}x_{2}\mathrm{d}x_{3}\mathrm{d}x_{4}\nonumber \\
 &=& stE\int_{-\infty}^{+\infty}\exp(-\lambda g(x,z))\exp(-(1-\lambda)g(y,z))\mathrm{d}\Phi(z),
\end{eqnarray}
which combining with (\ref{c231}) implies the desired result (\ref{2.10}).

(ii). Note that for $x\leq y$
\begin{eqnarray}
\label{c234}
&&P\left(\frac{1}{n}\widetilde{L}_{n}\leq s, \frac{1}{n}L_{n}\leq t, \widetilde{M}_{n}\leq u_{n}(x), M_{n}\leq u_{n}(y)\right)\nonumber\\
&&=P\left(\frac{1}{n}\widetilde{L}_{n}\leq s, \frac{1}{n}L_{n}\leq t, \widetilde{M}_{n}\leq u_{n}(x), \widehat{M}_{n}\leq u_{n}(y)\right)\nonumber\\
&&=P\left(\frac{1}{n}\widetilde{L}_{n}\leq s, \frac{1}{n}L_{n}\leq t, \widetilde{M}_{n}\leq u_{n}(x), \widehat{M}_{n}\leq u_{n}(y), \widetilde{M}_{n}\leq \widehat{M}_{n}\right)\nonumber\\
&&\ \ \ +P\left(\frac{1}{n}\widetilde{L}_{n}\leq s, \frac{1}{n}L_{n}\leq t, \widetilde{M}_{n}\leq u_{n}(x), \widehat{M}_{n}\leq u_{n}(y), \widetilde{M}_{n}> \widehat{M}_{n}\right)\nonumber\\
&&=P\left(\frac{1}{n}\widetilde{L}_{n}\leq s, \frac{1}{n}\widehat{L}_{n}\leq t, \widetilde{M}_{n}\leq u_{n}(x), \widehat{M}_{n}\leq u_{n}(y), \widetilde{M}_{n}\leq \widehat{M}_{n}\right)\nonumber\\
&&\ \ \ +P\left(\frac{1}{n}\widetilde{L}_{n}\leq \min\{s,t\}, \widetilde{M}_{n}\leq u_{n}(x), \widehat{M}_{n}\leq u_{n}(y), \widetilde{M}_{n}> \widehat{M}_{n}\right).
\end{eqnarray}
It follows from case (i) that for any $x,y\in \mathbb{R}$
\begin{eqnarray}
\label{c235}
&&\lim_{n\rightarrow\infty}P\left(\frac{1}{n}\widetilde{L}_{n}\leq s, \frac{1}{n}\widehat{L}_{n}\leq t, \widetilde{M}_{n}\leq u_{n}(x), \widehat{M}_{n}\leq u_{n}(y)\right)\nonumber\\
 &&= stE\int_{-\infty}^{+\infty}\exp(-\lambda g(x,z))\exp(-(1-\lambda)g(y,z))\mathrm{d}\Phi(z),
\end{eqnarray}
which indicates that $\widetilde{L}_{n}, \widehat{L}_{n}$ are asymptotically independent of $\widetilde{M}_{n}$ and $\widehat{M}_{n}$.
Thus, we have for $x<y$
\begin{eqnarray}
\label{c236}
&&\lim_{n\rightarrow\infty}P\left(\frac{1}{n}\widetilde{L}_{n}\leq s, \frac{1}{n}\widehat{L}_{n}\leq t, \widetilde{M}_{n}\leq u_{n}(x), \widehat{M}_{n}\leq u_{n}(y), \widetilde{M}_{n}\leq \widehat{M}_{n}\right)\nonumber\\
&&= \lim_{n\rightarrow\infty}P\left(\frac{1}{n}\widetilde{L}_{n}\leq s, \frac{1}{n}\widehat{L}_{n}\leq t\right)\times \lim_{n\rightarrow\infty} P\left(\widetilde{M}_{n}\leq u_{n}(x), \widehat{M}_{n}\leq u_{n}(y), \widetilde{M}_{n}\leq \widehat{M}_{n}\right)\nonumber\\
 &&= stE\int_{-\infty}^{+\infty}\exp(-\lambda g(x,z))\exp(-(1-\lambda)g(y,z))-\lambda\exp(-g(x,z))\mathrm{d}\Phi(z)
\end{eqnarray}
and
\begin{eqnarray}
\label{c237}
&&\lim_{n\rightarrow\infty}P\left(\frac{1}{n}\widetilde{L}_{n}\leq \min\{s,t\}, \widetilde{M}_{n}\leq u_{n}(x), \widehat{M}_{n}\leq u_{n}(y), \widetilde{M}_{n}> \widehat{M}_{n}\right)\nonumber\\
&&= \lim_{n\rightarrow\infty}P\left(\frac{1}{n}\widetilde{L}_{n}\leq \min\{s,t\}\right)\times \lim_{n\rightarrow\infty} P\left(\widetilde{M}_{n}\leq u_{n}(x), \widehat{M}_{n}\leq u_{n}(y), \widetilde{M}_{n}> \widehat{M}_{n}\right)\nonumber\\
 &&= \min\{s,t\}E\int_{-\infty}^{+\infty}\lambda \exp(-g(x,z))\mathrm{d}\Phi(z).
\end{eqnarray}
Plugging (\ref{c236}) and (\ref{c237}) into (\ref{c234}), we get the desired result (\ref{2.11}).

\bigskip



\end{document}